\documentclass[12pt]{article}

\usepackage[a4paper]{geometry}
\usepackage{ajc-FR}	
\usepackage{amsmath,amsfonts,latexsym}
\usepackage{graphicx}
\usepackage{hyperref}
\usepackage{placeins}

\Volume{XY(Z)}
\Year{2021}
\firstpage{1}
\lastpage{4}
\Revised{27 June 2021}
\runninghead{\copyright \, F. Rowley 2021\,/\,S(7)1696 DRAFT arXiv version, pages 1 -- 4}
\numberwithin{equation}{section}

\begin{document}

\title{\bf An Improved Lower Bound for $S(7)$\\ and Some Interesting Templates}
\author{Fred Rowley\thanks{formerly of Lincoln College, Oxford, UK.}}
\Addr{West Pennant Hills, \\NSW, Australia.\\{\tt fred.rowley@ozemail.com.au}}


\date{\dateline{27 June 2021}{DD Mmm CCYY}}

\maketitle


\begin{abstract}

In this simple paper, we exhibit a Schur partition giving rise to a triangle-free linear colouring of $K_{1697}$ in 7 colours.  Thus we show that the Schur number $S(7) \ge 1696$ and the multicolour Ramsey number $R_{7}(3) \ge 1698$.  

We also demonstrate a specific partition of the closed integer interval $[1, 376]$ which, through repeated use of a specific construction, allows us to conclude that $S(k+5) \ge 376.S(k)+160$ for any positive integer $k$.  The existence of this partition, with its specific properties, implies that $\lim_{\substack{r \rightarrow \infty}} R{_r}(3)^{1/r} > 3.273$.

\bigskip\noindent
\end{abstract}
\bigskip
\small DRAFT \copyright Fred Rowley  June 2021.
\bigskip

\section{Introduction}	

The first part of this paper arose from attempts to produce a general means of extending triangle-free linear Ramsey graphs by including one or more additional colours, in a way that reliably produced superior lower bounds for the relevant Ramsey numbers and  Schur numbers.  At this point, we note the isomorphisms between the 'distance sets' defining colourings of linear Ramsey graphs, and the corresponding Schur partitions.  

In this endeavour, it was noted that few of the most significant lower bounds for 'small' Ramsey graphs had been produced by combining graphs or partitions, with the notable exception of that given by Chung \cite{Chung} demonstrating that $R_{4}(3) \ge 51$.  However, Chung's construction does not depend on or preserve linearity in the colouring, so it is not useful in addressing the Schur numbers.  Neither the pioneering construction of Schur \cite{Schur} , nor the stronger result of Abbott and Hanson (Corollary 2.1 in \cite{AbbH}) produce strong lower bounds on $S(k)$ for $5 \le k \le 7$, leaving such bounds to be set by computer search programs, a very important example being those mentioned in \cite{FrSw}.  

This paper describes a search strategy that was used to find an improved upper bound for $S(7)$, and shows one particularly symmetrical result.  The partition can easily be verified as being sum-free, so no complex proof is needed.  The search was effective in a limited time, only because the strategy was constrained to a certain very limited range of patterns.   

For the second result, a proof in terms of linear graphs is already given in \cite{FR-GLRGC}.  The numerical result depends on the existence of a specific graph, for which the distance set (a Schur partion) is attached to this paper.  The implications of the existence of this partition for Schur numbers are briefly mentioned here.  

Some notation is defined in Section 2.  

In Section 3, we describe the search strategy that resulted in the attached partition, and the motivation for it. 

In Section 4, we describe a key partition of order 376, and show how it leads to the formula $S(k+5) \ge 376.S(k)+160$.

Section 5 contains some very brief remarks on further implications of these results.

\section{Definitions and Notation}	

In this paper:

A set of integers of the form $\{i \mid a \le i \le b\}$ for integers a, b, may be written as $[a, b]$.  

If a set of integers contains no three integers $a, b, c$, not necessarily distinct, such that $a + b = c$, then it is {\it sum-free}.  If a set $S$ is partitioned into disjoint non-empty subsets $S_i$, for $1 \le i \le k$, and all the subsets in the partition are sum-free, that partition is a {\it Schur partition}.  For integer $k > 0$, the {\it Schur number} $S(k)$ is the maximum value of $n$ such that a Schur partition of $[1, n]$ into $k$ subsets exists: its existence is established by Ramsey's Theorem. 

\section{A Simple Search Strategy}	

The highest existing lower bounds for $S(6)$ and $S(7)$ derive from work by Fredricksen and Sweet, whose search strategy is set out in \cite{FrSw}.  It can be paraphrased briefly as a process allocating numbers to subsets in a candidate partition, each allocation having regard to the number of 'blockages' (in their terminology) already applying to that number as a result of previous allocations.  Some optimisation was used to decide the order in which numbers were allocated, and many trials were carried out.  

Importantly for run-times, the authors used a constraint that the partition should be symmetrical (with special provision in the $S(6)$ case to allow the numbers 179 and 358 to be in different sets).  There was a further more flexible constraint (again paraphrasing) that in constructing candidates for $S(k)$, one subset was restricted to contain only numbers above a certain minimum, and their complements.  The chosen minimum was intended to be below, but close to, a known lower bound for $S(k-1)$.  

The partition in the file attached was produced by a more straightforward tree search, although the constraints in this case were much tighter.

The attempt to produce a higher lower bound for $S(7)$ was initiated partly because the value of 1680 was seen to be very close to the maximum of 1664 obtainable from the new construction described in \cite{FR-GLRGC}.  Since the value of 536 as a lower bound for $S(6)$ was well clear of the best attainable by compounding other graphs, $S(7)$ seemed a more natural target for improvement.  It was also noted that the ratio of 1680/536 falls well short of the ratio of 536/160 - noting that 160 is reported to be the precise value of $S(5)$.  

The strategy was to first decide on the order of partitions to be tested (in this case 1696) and then to take examples of an almost-symmetrical partition of a base set (in this case, $[1, 536]$), into six colours, putting the next integer (537) in a new seventh subset.  For each integer $x \le 537$ in the partition, its complement $(1697-x)$ was included in the same subset.  It was then only necessary to carry out a plain tree-search for satisfactory partitions, starting from 538 and working upwards.  It should be noted that the integers immediately above 537 are already heavily 'blocked' by previous entries in this structure. The fact that the starting partition of $[1, 536]$ is markedly non-critical appeared to help the search by limiting the number of blockages slightly.  

In this case, the resulting partition may be fully symmetrical, since 1697 is not divisible by 3, and the example provided in the attached ancillary file was chosen to be fully symmetrical.  However, we note that a profusion of examples were found that were not symmetrical, and the attached partition is highly non-critical: which may be a sign that there are larger partitions to be found.

\section{Three Particular Partitions}	

This section briefly describes the attributes of a sum-free partition of $[1, 376]$ into 6 subsets, which -- if converted into a linear Ramsey graph, and borrowing the terminology of \cite{FR-GLRGC} -- has a 'triangle-free template' and a parameter $\phi = 160$.  Its profile is included in the ancillary file.  

The application of the Ramsey graph result to Schur partitions, mentioned in that paper, is a simple consequence of an effective equivalence between the so-called 'linear' graphs and Schur partitions -- the distance sets of the linear graph being isomorphic to the subsets comprising the Schur partition, and similarly sum-free.  

This partition has been verified as suitable for use in the construction described in \cite{FR-GLRGC}, which assures that the constructed partition remains sum-free.  Used with the strong Schur partition of $[1, 4]$ into subsets $\{1, 4\} , \{2, 3\}$, it produces the partition of $[1, 1664]$ described above, and attached in the ancillary file.  It does this by replacing the subset equivalent to the 'template colour' with multiple translations of that subset following the pattern of that second partition, as also described in \cite{FR-GLRGC}.

In fact, this partition can be combined with any strong Schur partition in this way, and therefore by repeated application, its existence implies that $S(k+5) \ge 376.S(k)+160$ for all integer $k > 0$.  By constructing an infinite series of such graphs, we can demonstrate that $$\lim_{\substack{r \rightarrow \infty}} R{_r}(3)^{1/r} \ge \sqrt[5]{376} = 3.273 \dots $$. 

This is a notable improvement over previous lower bounds, of which the previous highest value of $3.199 \dots$ is to be found in \cite{XXER}.

A sum-free partition of $[1, 109]$ into 5 subsets with parameter $\phi = 39$, constructed similarly, is also included in the ancillary file, and by the same logic, its existence implies that $S(k+4) \ge 109.S(k)+39$ for all integer $k > 0$.  

A sum-free partition of $[1, 33]$ with similar properties in 4 colours, and with parameter $\phi = 6$, is included in the same ancillary file.  Its existence similarly implies that $S(k+3) \ge 33.S(k)+6$ for all integer $k > 0$.  

Each of these partitions can of course be applied as described in \cite{Rowley3} to produce weak Schur partitions of unlimited size.  


\section{Some Remarks}	
	
As regards Section 3, it is interesting to note that many satisfactory partitions of order 1696 were found relatively quickly -- in fact, on the first day of test-running of the search programs.  However, no larger partition into 7 colours has been found so far, despite further searches for sets of order close to 1696.  It is interesting to note that, like 1680, the integer 1696 is divisible by 16; which perhaps may be a sign of some structure that has not been detected in either of the partitions.  

As regards Section 4, the main objective in the current paper is to make clearer the structure of certain graphs used to produce previous results, by providing a set of examples.  Some of the partitions exhibited here are equivalent to the linear graph colourings underpinning the results in \cite{FR-GLRGC}, as well as some of the results mentioned in \cite{Rowley3}. It is hoped that the availability of these examples will be useful to other researchers in this area.






\end{document}